\documentclass{article}

\usepackage{color, amsmath, amssymb, amsthm}
\usepackage[a4paper,body={16cm,23cm}]{geometry}

\definecolor{dblue}{rgb}{0,0,0.7}
\newtheoremstyle{mythm}{11pt}{11pt}{\it\color{dblue}}{}{\bf\color{dblue}}{.}{ }{}
\theoremstyle{mythm}
\newtheorem{thm}{Theorem}

\begin{document}
\title{The sub-leading coefficient of the $L$-function of an elliptic curve}
\author{Christian Wuthrich}
\maketitle
\abstract{We show that there is a relation between the leading term at $s=1$ of an $L$-function of an elliptic curve defined over an number field and the term that follows.}
\medskip


%
 %

Let $E$ be an elliptic curve defined over a number field $K$. We will assume that the $L$-function $L(E,s)$ admits an analytic continuation to $s=1$ and that it satisfies the functional equation. By modularity~\cite{modularity}, we know that this holds when $K=\mathbb{Q}$. The conjecture of Birch and Swinnerton-Dyer predicts that the behaviour at $s=1$ is linked to arithmetic information. More precisely, if
\begin{equation*}
 L(E,s) = a_r \, (s-1)^r + a_{r+1}\, (s-1)^{r+1} + \cdots
\end{equation*}
is the Taylor expansion at $s=1$ with $a_r\neq 0$, then $r$ should be the rank of the Mordell-Weil group $E(K)$ and the leading term $a_r$ is equal to a precise formula involving the Tate-Shafarevich group of $E$. It seems to have passed unnoticed that the sub-leading coefficient $a_{r+1}$ is also determined by the following formula.

\begin{thm}
 With the above assumption, we have the equality
 \begin{equation}
  a_{r+1} = \Bigl( [K:\mathbb{Q}]\cdot(\gamma + \log(2\pi)) - \tfrac{1}{2}\log(N) - \log \vert\Delta_K\vert \Bigr) \cdot a_r
 \end{equation}
 where $\gamma=0.577216\dots$ is Euler's constant, $N$ is the absolute norm of the conductor ideal of $E/K$ and $\Delta_K$ is the absolute discriminant of $K/\mathbb{Q}$.
\end{thm}

In particular, the conjecture of Birch and Swinnerton-Dyer also predicts completely what the sub-leading coefficient $a_{r+1}$ should be. One consequence for $K=\mathbb{Q}$ is that for all curves with conductor $N>125$, and this is all but $404$ isomorphism classes of curves, the sign of $a_{r+1}$ is the opposite of $a_r$. Of course, it is believed that $a_r$ is positive for all $E/\mathbb{Q}$.

\begin{proof}
 Set $f(s) = B^s \cdot \Gamma(s)^n$ with $n=[K:\mathbb{Q}]$ and $B=\sqrt{N} \cdot \vert \Delta_K\vert /(2\pi)^n$. Then $\Lambda(s) = f(s) \cdot L(E,s)$ is the completed $L$-function, which satisfies the functional equation $\Lambda(s) = (-1)^r\cdot \Lambda(2-s)$, see~\cite{husemoller}. For $i\equiv r+1\pmod{2}$ it follows that $\frac{d^i}{ds^i}  \Lambda(s)\bigr\vert_{s=1} = 0$. Hence for $i=r+1$, we obtain that
 \begin{equation*}
  (r+1)\cdot f'(s) \cdot \frac{d^r}{ds^r}L(E,s) + f(s)\cdot \frac{d^{r+1}}{ds^{r+1}} L(E,s)
 \end{equation*}
 is zero at $s=1$. Therefore $(r+1) \, f'(1) \,r!\, a_r + f(1)\,(r+1)!\, a_{r+1} = 0$. It remains to note that $f(1) = B$ and $f'(1) = B\cdot \bigl( \log(B) + n\cdot \Gamma'(1) \bigr)$ together with $\Gamma'(1) = - \gamma$.
\end{proof}

Obviously a similar formula holds for the $L$-function of a modular form of weight $2$ for $\Gamma_0(N)$. More generally, for any $L$-function with a functional equation there is a relation between the leading and the sub-leading coefficient of the Taylor expansion of the $L$-function at the central point.

Sub-leading coefficients of Dirichlet $L$-functions have been investigated; for instance Colmez~\cite{colmez} makes a conjecture, which is partially known. However these concern the much harder case when $s$ is not at the centre but the boundary of the critical strip of the $L$-function.


\bibliographystyle{amsplain}
\bibliography{sl}

\end{document}